\documentclass[11pt]{amsart}
\usepackage{graphicx}
\usepackage{epstopdf}
\font\caps=cmcsc10
\newcommand{\br}{\ensuremath{\mathbb R}}
\newcommand{\fig}[2]{\includegraphics[scale=#2]{#1.eps}}
\newcommand{\figno}[1]{\smallskip\centerline{Figure #1}\medskip}
\newcommand{\place}[3]
 {\text{\kern#1pt
    \smash{\raise#2pt\hbox{\rlap{#3}}}
  \kern-#1pt\kern-1ex}}

\begin{document}

\title
{A topological menagerie}
\author
{Paul Melvin}
\date{}
\address{Bryn Mawr College, Bryn Mawr, PA 19010-2899}
\email{pmelvin@brynmawr.edu}
\maketitle

Some of the party games that young teenagers play have surprisingly rich mathematical content.  ``Entanglement" is one such game in which couples are linked with strings tied tightly to their wrists, as shown in Figure 1(a)

\medskip
\centerline{\fig{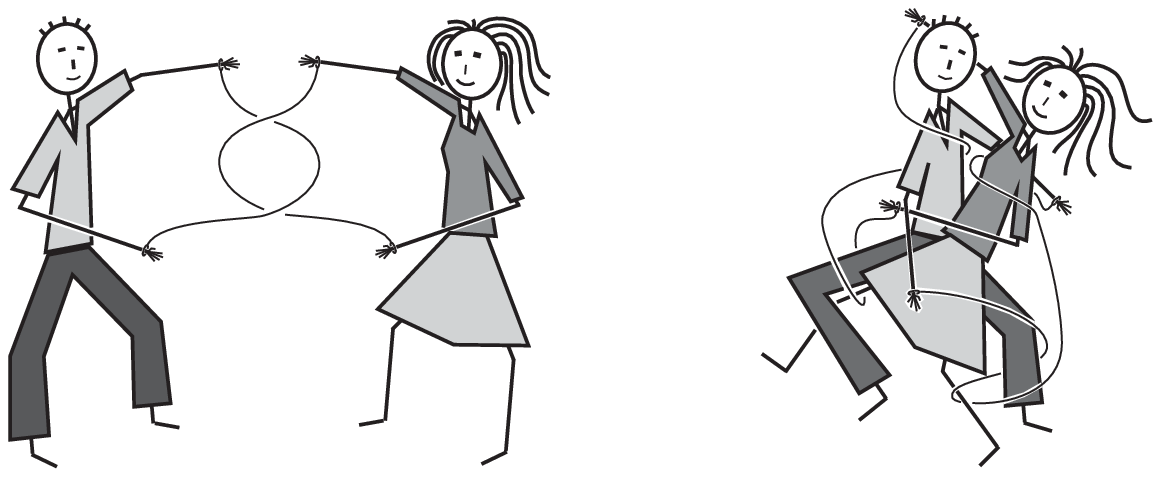}{.75}}
\centerline{\kern 10 pt(a) \kern 120pt (b)}

\figno{1: Entanglement}

\noindent
and then challenged to disentangle themselves, often leading to Figure 1(b).  In fact there is a quick solution (needless to say missing the point of the game): simply push a bit of the girl's string under the boy's at his wrist, pass the resulting loop over his hand, and then pull it free on the other side, as shown in Figure 2.

\medskip
\centerline{\fig{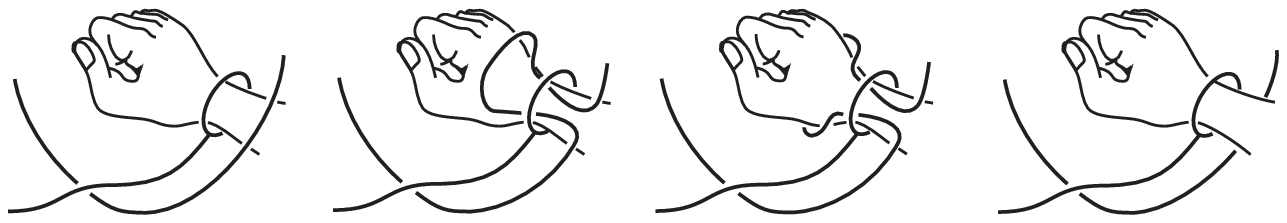}{.8}}

\figno{2: Disentanglement}

One way to formulate this mathematically is to ask for an isotopy to remove the meridian loop $m$ from a rigid wire $w$, embedded in the $3$-sphere as shown in Figure 3(a).  (The wire represents the boy with his attached string, where the upper loop of the wire is his right hand, while the meridian represents the girl with her string.)  Of course this is easy: just slide $m$ along $w$ to a point just below the upper loop of $w$, and then stretch it out and pull it off the top.  Or put differently, observe that $w$ is isotopic to the trivial wire $\circ$\hskip-.025in$-\!\!\!-\!\!\!-\!\!\!-$\hskip-.02in$\circ$ and so the meridian slips right off.

\medskip
\centerline{\fig{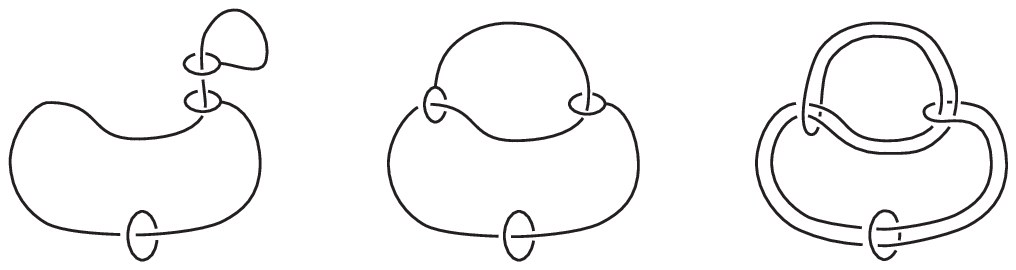}{1}}
\place{57}{52}{$w$}
\place{56}{5}{$m$}
\place{162}{52}{$w$}
\place{159}{5}{$m$}
\place{263}{52}{$k$}
\place{259}{5}{$m$}

\centerline{(a) Entanglement  \kern 30pt (b) Coffin's puzzle \kern 20pt (c) The knot $k$ }

\figno{3}

A similar puzzle attributed to Steward Coffin \cite{db} was considered by Inta Bertuccioni in a recent issue of the Monthly \cite{bertucionni}.  Here the wire $w$ is configured slightly differently (as in Figure 3(b)) and there is in fact no solution, that is $m$ cannot be isotoped off $w$.  Bertuccioni proves this by an explicit calculation showing that $m$ represents a non-trivial element in the fundamental group of the complement of $w$.

The purpose of this note is to give a more conceptual proof of the impossibility of solving Coffin's puzzle, and then to generalize this proof, using an elementary but non-trivial result in knot theory, to show that all but one of the vast ``menagerie" of possible puzzles suggested by entanglement and Coffin's puzzle are also unsolvable.

First consider Coffin's puzzle.  Note that the wire $w$ consists of two unknots $\lambda$ joined by an arc $\alpha$.  From this one obtains a knot $k$ by banding the components of $\lambda$ together along $\alpha$, as shown in Figure 3(c).  Actually there are infinitely many knots that can be formed in this way as there may be twists in the band; the one pictured is readily seen to be the square knot.  

Now observe that if $m$ could be isotoped off $w$, then it certainly could be isotoped off $k$, since the complement of $k$ contains the complement of $w$ (with $\alpha$ thickened a bit).  But then the lift $\tilde k$ of $k$ to the infinite cyclic cover of the complement of $m$ would consist of an infinite number of copies of $k$.  (Note that the complement of $m$ is a solid torus, and the cover is an infinite solid cylinder.)  However, it is easily seen that the components of $\tilde k$ are {\sl unknotted}, contradicting the fact that $k$ is a nontrivial knot (the square knot).  Indeed, viewing $k$ as the closure of a tangle with axis $m$, the link $\tilde k$ is obtained by composing infinitely many copies of this tangle, as shown in Figure 4 (cf.\ \cite{rolfsen}).    

\medskip
\centerline{\fig{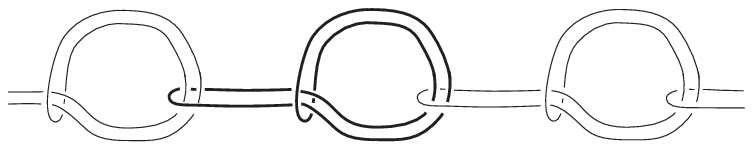}{1}}
\figno{4: $\tilde k$ in the infinite cyclic cover}

\noindent
Thus each component of $\tilde k$ looks like a long worm that can be shrunk by pushing from its ``free" end. 

\break 

For a general puzzle in the menagerie we allow any embedding of the wire for which $\lambda$ is an unlink; the arc $\alpha$ can be arbitrary.  One such puzzle is shown in Figure 5(a).  

\medskip
\centerline{\fig{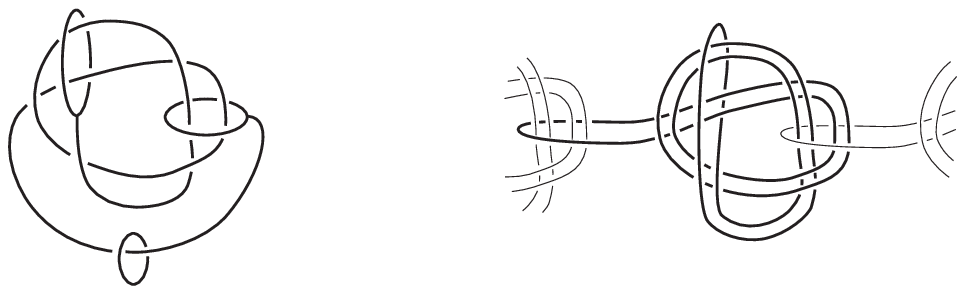}{1}}
\place{70}{52}{$w$}
\place{68}{5}{$m$}

\centerline{(a) A beast in the menagerie \kern 30pt (b) The associated lift $\tilde k$\kern 10pt}

\figno{5}

\noindent
The argument above shows that there is no hope for a solution unless the associated knot $k$ is trivial, since each component of the lift $\tilde k$ is unknotted (as seen in Figure 5(b); in general slide $m$ close to an endpoint of $\alpha$ before taking the cover).  However, by a theorem of Marty Scharlemann \cite{scharlemann}, $k$ is trivial if and only if the arc $\alpha$ is isotopic (fixing $\lambda$) to a trivial arc.  Thus ``most" of the puzzles in the menagerie, indeed all but those equivalent to entanglement,  have no solution.  

From a practical point of view one might ask how to recognize when $\alpha$ is trivial.  This is actually a difficult question.  However there is a simple test for the triviality of the {\sl homotopy} class of $\alpha$:  The fundamental group of the complement of $\lambda$ is free of rank two, with generators $x$ and $y$ corresponding to the two components of $\lambda$, and $\alpha$ (oriented from the $x$-loop to the $y$-loop of $\lambda$) determines a word in $x$ and $y$.  Since $\alpha$ can be made to spiral around $\lambda$ at it's endpoints, its homotopy class corresponds in fact to an {\sl equivalence class} of words, where $u\sim v \iff u=x^mvy^n$ for some integers $m,n$.  That is {\sl $\alpha$ is homotopically trivial if and only if $($the word$)\ \alpha\sim 1$}.  For example, the entanglement puzzle has $\alpha = xy \sim 1$ as expected, while Coffin's puzzle has $\alpha = yx \not\sim 1$, and the puzzle in Figure 5 has $\alpha = yxyx \not\sim 1$.  

Of course there are many homotopically trivial arcs $\alpha$ that are isotopically nontrivial, or equivalently whose associated knots are nontrivial.  Consider for example the  puzzle in Figure 6, which has $\alpha = xy \sim 1$. 

\medskip
\centerline{\fig{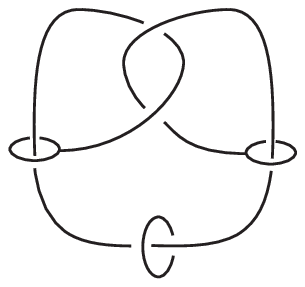}{1}}

\figno{6: A homotopically trivial beast}

\noindent
The associated knot has Jones polynomial \cite{jones} $t^{-5}-t^{-4}-t^{-1}+2-t+t^2+t^5-t^6$ (the Alexander polynomial is trivial) and so this puzzle has no solution.

\end{document}